\documentclass[letterpaper, journal,onecolumn, draftcls]{IEEEtran}  

\IEEEoverridecommandlockouts                              
\usepackage{url}
\usepackage{graphicx}
\usepackage{graphics} 
\usepackage{epsfig} 
\usepackage{mathptmx} 
 \usepackage{mathrsfs}
\usepackage{times} 
\usepackage{amsmath} 
\usepackage{amssymb}  
\usepackage{color}
\usepackage{float}
\usepackage{algorithm}
\usepackage{algorithmic}
\usepackage{etoolbox}

\makeatletter
\def\hlinew#1{%
  \noalign{\ifnum0=`}\fi\hrule \@height #1 \futurelet
   \reserved@a\@xhline}
\patchcmd{\@makecaption}
{\scshape}
{}
{}
{}
\makeatother

\title{Forecasting Spatio-Temporal Renewable Scenarios: a Deep Generative Approach}

\author{Congmei Jiang, Yize Chen, Yongfang Mao, Yi Chai and Mingbiao Yu

\thanks{C. Jiang, Y. Mao, Y. Chai are with College of Automation, Chongqing University. Emails: jiangcongmei@gmail.com and \{yfm,chaiyi\}@cqu.edu.cn. Y. Chen is with Department of Electrical and Computer Engineering, University of Washington. Email: yizechen@uw.edu. M. Yu is with School of Instrument Science and Engineering, Southeast University. Email: ymb\_moon@126.com.}}

\begin{document}

\maketitle
\begin{abstract}
The operation and planning of large-scale power systems are becoming more challenging with the increasing penetration of stochastic renewable generation.  In order to minimize the decision risks in power systems with large amount of renewable resources, there is a growing need to model the short-term generation uncertainty.  By producing a group of possible future realizations for certain set of renewable generation plants, scenario approach has become one popular way for renewables uncertainty modeling. However, due to the complex spatial and temporal correlations underlying in renewable generations, traditional model-based approaches for forecasting future scenarios often require extensive knowledge, while fitted models are often hard to scale. To address such modeling burdens, we propose a learning-based, data-driven scenario forecasts method based on generative adversarial networks (GANs), which is a class of deep-learning generative algorithms used for modeling unknown distributions. We firstly utilize an improved GANs with convergence guarantees to learn the intrinsic patterns and model the unknown distributions of (multiple-site) renewable generation time-series. Then by solving an optimization problem, we are able to generate forecasted scenarios without any scenario number and forecasting horizon restrictions. Our method is totally model-free, and could forecast scenarios under different level of forecast uncertainties. Extensive numerical simulations using real-world data from NREL wind and solar integration datasets validate the performance of proposed method in forecasting both wind and solar power scenarios.

\end{abstract}


\section{Introduction}
\label{intro}
Decision-making under uncertainty environments have long been a challenging problem for power system engineers. There is no perfect forecasts, and thus consideration of generation uncertainties is a necessity for reliable, efficient power system operation planning~\cite{bertsimas2013adaptive}. To accommodate for the higher penetration of renewable energy such as wind, solar and hydro power, it becomes an urgent need for developing more accurate, scalable and efficient approaches to model short-term renewable generation uncertainties~\cite{zhu2015optimization}. Handling the uncertainty for renewable generations is key to increasing economic benefits and enforcing reliability criteria for decision-making. One widely used approach to capture the uncertainties in renewable resources is by modeling the future via generating a set of possible future time-series called \emph{scenarios}. Comparing to other uncertainty modeling techniques such as probabilistic forecasts and quantile forecasts, scenarios reflect the joint distributions of renewables generation at different locations and varying lead time~\cite{pinson2012evaluating}. Scenario approach has been playing an important role in a series of stochastic and robust optimization problems such as unit commitment, energy trading strategy, storage sizing and etc~\cite{wang2008security,ma2013scenario,wang2017look}.

Despite such promise and wide applications, generating scenarios that are able to accurately inform system operators on power generation uncertainties remains to be a challenging problem. One of the biggest challenges of scenario generation is the difficulty of modeling or learning the unknown  stochastic processes that drive renewable power generation.  Researchers have conducted extensive research on using probabilistic, statistical models for scenario generations. In~\cite{pinson2012evaluating, pinson2009probabilistic,tang2018efficient}, Gaussian copula or pair-copula model are used to generate statistical scenarios that accounts for both the interdependence structure of prediction errors and the predictive distributions from wind power probabilistic forecasting.  Auto-regressive moving average (ARMA) along with Monte Carlo simulation are used to generate wind power scenarios in~\cite{meibom2011stochastic}. For generating spatial correlated scenarios, time series models~\cite{morales2010methodology,le2015probabilistic} are illustrated to produce a set of plausible scenarios characterizing the uncertainty associated with wind speed at different geographic sites. 

However, most approaches mentioned above required large amount site-specific modeling knowledge on renewable resources. What's more, the model assumptions on future power generations, e.g. multi-dimensional Gaussian assumptions on forecast horizons, are normally not held and vary by locations and time. The intermittent and time-varying nature of renewables, the complex spatial and temporal interactions make most of these methods difficult to apply and hard to scale in practice. Thus generated scenarios may not represent the intrinsic patterns and realistic time-series of real historical observations of renewable energy resources.

Instead of modeling the stochastic processes explicitly using statistical models, data-driven methods have also been considered to represent the complex dynamics in renewables generation processes~\cite{liu2018small,qin2019optimized}. In~\cite{vagropoulos2016ann}, a scenario generation methodology based on artificial neural networks (ANNs) is proposed to firstly get an accurate point forecasts, yet such supervised approach heavily relies on the Gaussian assumption and statistical information of forecast errors.

On the other hand, generative models have been used to directly learn the underlying distribution of renewables generation processes. Generative models are trained to learn the true data distribution during the training phases. The learning process is unsupervised, and once the model is trained, it can efficiently generate data that has a similar distribution to the original ones~\cite{li2018novel}.  In~\cite{zhanga2018optimized}, variational autoencoder~(VAE) is adopted to generate scenarios for wind and PV power. 
In~\cite{chen2018model, jiang2018scenario}, the authors firstly propose to use generative adversarial networks~(GANs)~\cite{goodfellow2014generative}, which is a model-free, data-driven and scalable approach for generating renewable scenarios by deep generative models. However, such method cannot incorporate forecast information and be applied for future scenario generations. A followup work~\cite{chen2018unsupervised} proposed to generate a group of future realizations based on trained GANs, yet generating scenarios for multiple spatially correlated sites are not discussed. 

\begin{figure}[h]
	\centering
	\includegraphics[width=0.8 \columnwidth]{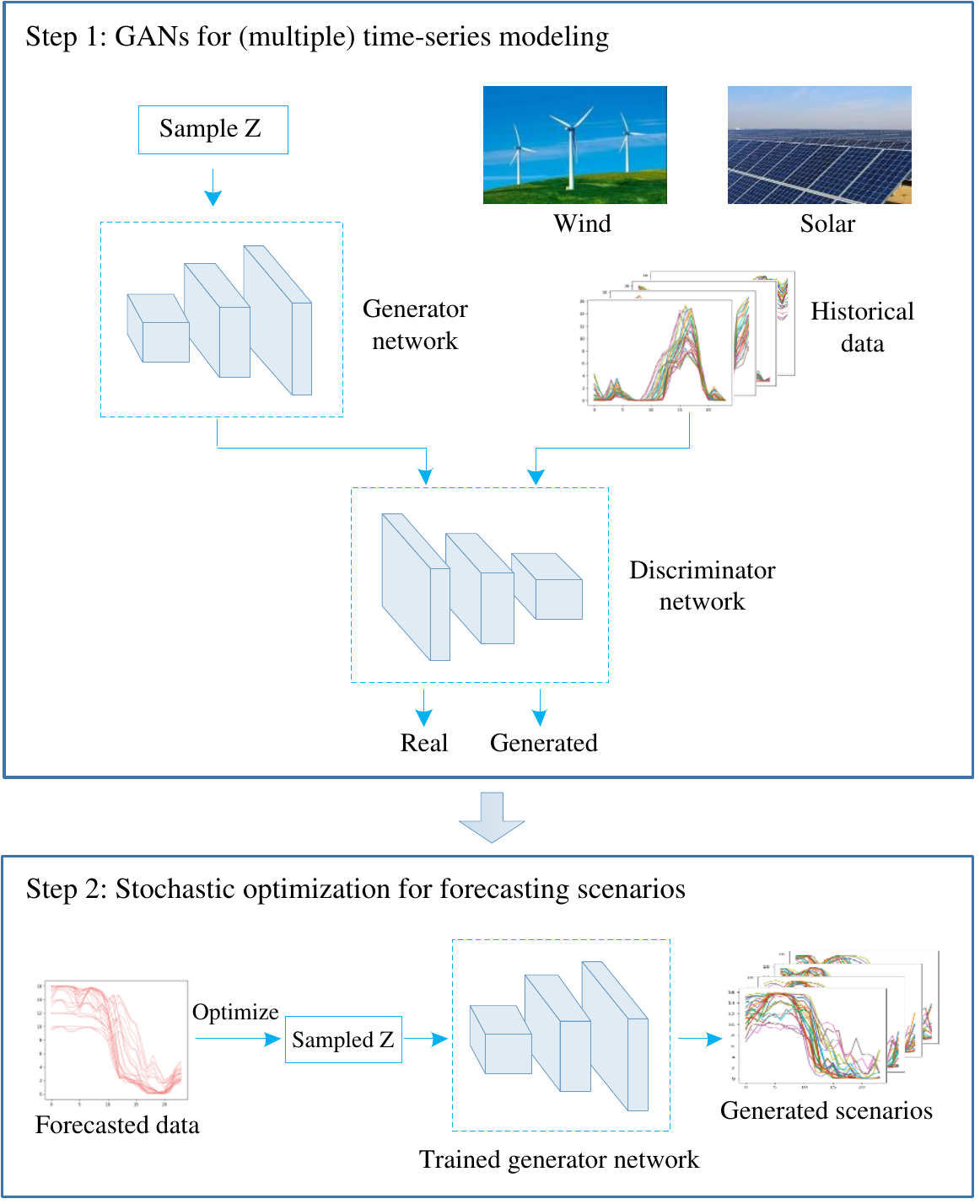}
	\caption{Illustration of the framework for different forecasting tasks based on GANs. In the first step, GANs are utilized to generate realistic scenarios for renewable generation. We then formulate an optimization problem in the second step based on point forecasts to generate future scenarios without number and time length restrictions.}
	\label{fig:schematic}
\end{figure}

Due to the correlations of meteorological conditions, renewable generation outputs of different locations exist certain unknown, hard-to-model correlations. The spatial dependence along with temporal correlations are both imperatives for joint uncertainty modeling, especially for power flow optimizations and transmission risk assessments. In this paper, we follow the line of deep generative model, and propose a novel method based on GANs to directly generate future scenarios for multiple-site renewable power generations. GANs are composed of two deep neural networks: a generator network trying to generate realistic samples, and a discriminator network trying to discriminate the input samples. In the training process, the generator network tries to ''fool" the discriminator network by generating realistic samples, while the discriminator network tries to distinguish the real training data from the output of the generator network. The joint training of this two networks form a minimax game. GANs lie in the category of unsupervised learning model, while trained generators can generate realistic scenarios resembling to renewable power generation samples. We observe the training instability and low generation quality problems, and propose techniques to improve GANs training. 

Our method for scenarios forecasts contains two steps. In the first step, we train GANs based on historical observations of renewable generations. Once training is completed, by formulating an optimization problem based on given point forecasts, we optimize over the noise vectors to find the future scenarios from trained generator's outputs. The only information we need is a historical dataset comprised of power generations of target sites, along with any available point forecasts. The proposed approach is free of any statistical assumptions, and can forecast scenarios without relying on any sampling techniques. Fig.~\ref{fig:schematic} illustrates the algorithmic framework for our proposed method.

Specifically, we make the following contributions:
\begin{itemize}
	\item Based on any provided point forecasts, our method is able to generate short-term forecasting scenarios with high-level flexibility on number of renewable generation sites and scenarios as well as the length of forecast horizons. To our knowledge, this is the first work that utilizes deep generative models for forecasting spatio-temporal scenarios.
	\item The improved generative model avoid the problems of exploding or vanishing gradients during training and can also achieve faster convergence to reduce the training time for the generative models.
	\item  Scenarios forecasted by our method can not only represent a group of future realizations, but also reflect the intrinsic complex spatial and temporal patterns lying in the renewable generation processes. Extensive numerical simulations validate the performance of proposed method on varying scenario forecasts tasks. 
\end{itemize}

The rest of this paper is organized as follows. Section~\ref{gan} describes the model setup for GANs along with improved training techniques. In Section~\ref{forecast} we describe our proposed method for using GANs to forecast spatiotemporal scenarios. Section~\ref{network} provides the model structure and training algorithms of GANs for renewable scenario generation. In section~\ref{experiment}, numerical simulations are conducted to validate the proposed technique for forecasting wind or solar scenarios through a comprehensive analysis comprising for both single site and spatial-correlated multiple sites. Conclusion marks are made in Section~\ref{conclusion}.

\section{Data-Driven Generative Model}
\label{gan}
In this section, we describe the formulation and training techniques for GANs, and illustrate how to utilize GANs as an efficient module for scenario generation. Improved training techniques for GANs, including using Wasserstein distance as loss metric, enforcing Lipschitz constraints and continuity conditions are discussed.


\subsection{Wasserstein GAN (WGAN)}
GANs utilizes the power of deep learning to learn the unknown target distribution. It can be regarded as a two-player zero-sum game between the two interconnected neural networks, the generator $G$ and the discriminator $D$, under the adversarial learning framework. Given a training set of historical renewable generation data, the generator's goal is to find a function that transforms a sample from known noise distribution to a sample following the same distribution as the historical observations. The discriminator's goal is to distinguish whether the input data comes from the generator or comes from real historical samples. When the adversarial networks are trained to an equilibrium, the discriminator can no longer distinguish between generated and historical data, meaning the generator can produce realistic samples as if they are coming from the true distribution.

Suppose the distribution of the historical data $X$ is represented by the probability density function $\mathbb{P}_r$, and a noise vector $Z$ is sampled from a given Gaussian distribution $\mathbb{P}_Z$ with known mean and variance\footnote{In training implementations, any known, easy-to-sample distribution can be used for GANs training.}. Just like any training procedure for neural networks model, we need to define the loss functions to guide the parameter updates.

We firstly formulate the loss function $L_G$ to update the weights of $G$'s parameters. During training, a batch of samples drawn with distribution $\mathbb{P}_Z$ output newly generated data under generated samples' distribution $\mathbb{P}_G$. A small $L_G$ can be achieved by maximizing $D(G(Z))$, which indicates the generated samples from distribution $\mathbb{P}_G$ look like real samples from discriminator's perspective. Following this guideline, the loss function is defined as
\begin{equation}
L_G=-\mathbb{E}_{Z\sim \mathbb{P}_Z}[D(G(Z))]
\end{equation}

The discriminator takes input samples either coming from generator or from real historical data. It is alternately trained with the generator. During training, the discriminator's goal is to distinguish between $\mathbb{P}_X$ and $\mathbb{P}_G$. In other words, to maximize the value between  $\mathbb{E}[D(\cdot)]$ and $\mathbb{E}[D(G(\cdot))]$. To update the parameters of $D$, the loss function $L_D$ can be similarly defined by
\begin{equation}
L_D=-\mathbb{E}_{x\sim\mathbb{P}_r}[D(X)]+\mathbb{E}_{Z\sim\mathbb{P}_Z}[D(G(Z))]
\end{equation}

As for the adversarial training of the two interconnected neural networks, the discriminator outputs a continuous value to measure the input samples. For a given $D$, maximized output $G(Z)$ means to minimize $-\mathbb{E}[D(G(Z))]$, while the discriminator wants to minimize $\mathbb{E}[D(G(Z)]$ for generated samples and maximize $\mathbb{E}[D(X)]$ for real samples. With the two loss functions $L_D$ and $L_G$ defined, we then can formulate the two-player game with a value function $V(G, D)$:

\begin{equation}
\label{equ:minmax}
\min_G\max_D \; V(G,D)=\mathbb{E}_{X\sim\mathbb{P}_r}[D(X)]-\mathbb{E}_{Z\sim \mathbb{P}_Z}[D(G(Z))]
\end{equation}

The minimax objective~\eqref{equ:minmax} can be theoretically interpreted as the dual of the Wasserstein distance for $\mathbb{P}_r$ and $\mathbb{P}_G$. The original form of Wasserstein distance is defined as follows:
\begin{equation}
\label{equ:wass}
W(\mathbb{P}_r,\mathbb{P}_G)=\inf_{\psi\in \prod(\mathbb{P}_r, \mathbb{P}_G)}\; \mathbb{E}_{(X,Y)\sim\psi}||X-Y||
\end{equation}
where $ \prod(\mathbb{P}_r, \mathbb{P}_G)$ denotes the set of all joint distributions $\psi(X,Y)$ whose marginals are $\mathbb{P}_r$ and $\mathbb{P}_G$, and $inf$ requires to find the joint distribution that have smallest distance for $X$ and $Y$.  Yet directly minimizing \eqref{equ:wass} is impractical for the generated samples parameterized by neural networks $G$. By borrowing Kantorovich-Rubinstein duality~\cite{villani2008optimal}, we could alternatively optimize over the dual form of the Wasserstein distance:
\begin{equation}
\label{equ:dual}
W(\mathbb{P}_r,\mathbb{P}_G)=\sup_{D}\; \mathbb{E}_{X\sim\mathbb{P}_r}[D(X)]-\mathbb{E}_{Z\sim \mathbb{P}_Z}[D(G(Z))]
\end{equation}

Wasserstein distance is also known as the Earth-Mover (EM) distance~\cite{villani2008optimal}, which indicates how much ``mass" is needed to transport from one distribution to the targeted distribution. In terms of model training, this distance has nice properties on indicating the distance between generated samples and real samples compared with other standard distance measurements~(e.g., Jensen-Shannon divergence, Kullback-Leibler divergence). We will show training convergence results and high-quality generated renewable scenarios in latter sections.

\begin{figure}[h]
	\centering
	\includegraphics[width=0.8 \columnwidth]{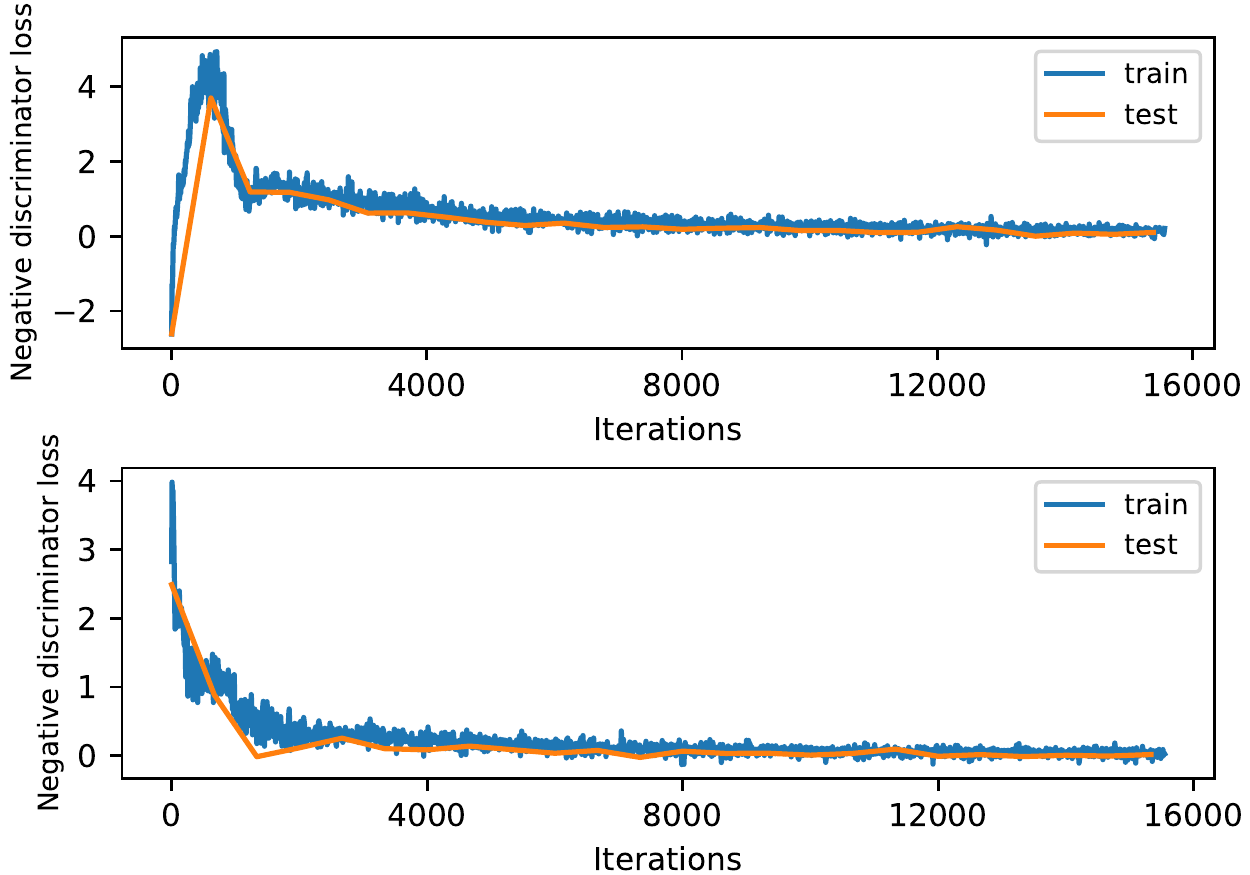}
	\caption{Training evolutions(blue) for GANs on (a) a wind dataset and (b) a solar dataset, respectively. The discriminator's losses are checked in every 100 generator iterations.}
	\label{fig:loss}
\end{figure}

\subsection{Improved Training Techniques}
WGAN uses Wasserstein distance to measure the distance difference between generated scenarios and historical data distribution, which theoretically addresses the problem of training convergence. The following techniques improve the constraints upon model weights, and guarantee quicker convergence and better scenario generation performances.

\subsubsection{Enforcing Lipschitz Constrains}
We follow the improved strategy proposed for imposing the Lipschitz constraint~\cite{gulrajani2017improved}. Inspired by
the optimal discriminator that has unit gradient norm almost everywhere under $\mathbb{P}_G$ and $\mathbb{P}_r$, the gradient penalty is given by
\begin{equation}
\label{equ:penalty}
GP(\tilde{X})=\mathbb{E}[(||\nabla_{\tilde{X}}D(\tilde{X})||_2-1)^2]
\end{equation}

where $\tilde{X}=tx +(1-t)G(Z)$ for $t\sim U[0,1]$. Given that enforcing the unit gradient norm constraint everywhere is intractable, this alternative way is an effective way to use for model training.

The gradient penalty term GP performs better than the standard weight clipping for Lipschtiz
constraint. The modified loss function stabilizes the GANs training over a wide range of architectures with almost no hyper-parameter tuning and can generate higher quality samples on different datasets~\cite{gulrajani2017improved}.

\subsubsection{Enforcing Continuity Conditions}
Since the gradient term can only be punished at sampled data points in the training process, a large part of the data points will not be sampled at all. In addition, the output of the generator is significantly different from the actual data point at the start of the training. The 1-Lipschtiz constraint is not enforced until the data distributions $\mathbb{P}_G$ and $\mathbb{P}_r$ are close enough to each other. To overcome these issues, an additional consistency term (CT) is proposed to improve the training~\cite{wei2018improving}. Instead of focusing on particular data points sampled on specific data points, a region around the real data manifold is considered. In particular, two perturbed data points $X'$  and $X''$ near observed real data point $X$ are used to check the continuity condition. The discriminator is Lipschitz continuous if there exists $M\geq 0$, such that for any input $X'$, $X''$, we must have
\begin{equation}
\label{equ:consistent}
||D(X')-D(X'')||_2\leq M ||X'-X''||_2
\end{equation}
where $X'$  and $X''$ can be found by applying the stochastic dropout to the hidden layers of the discriminator. GANs performance can be further improved by controlling the output of second-to-last layer of the discriminator $D^{(-2)}(\cdot)$. The final consistency regularization takes the following form to penalize the violation of \eqref{equ:consistent}:
\begin{equation}
\label{equ:consistency}
\begin{split}
CT(X', X'')=\mathbb{E}_{X\sim \mathbb{P}_r}[&\max (0, \;
||D^{(-2)}(X')-D^{(-2)}(X'')||_2+\\&0.1||D(X') - D(X'')||_2-M')]
\end{split}
\end{equation}
where $M'$ is a bounded constant.

In summary, the gradient penalty term GP~\eqref{equ:penalty} enforces the continuity over the points sampled between the real and generated points; while the consistency term CT~\eqref{equ:consistency} can complement the former by focusing on the region around the real data manifold instead. Therefore, these two terms can be used together to improve the training of GANs. Putting all together, the improved objective function for WGAN training can be expressed as
\begin{equation}
\label{equ:new_GAN}
\begin{split}
\min_G\max_D \; V(G,D)=\mathbb{E}_{X\sim\mathbb{P}_r}[D(X)]-\mathbb{E}_{Z\sim \mathbb{P}_Z}[D(G(Z))]\\ - \lambda_1 GP(\tilde{X})-\lambda_2 CT(X', X'')
\end{split}
\end{equation}
where $\lambda_1$ and $\lambda_2$ are used to balance the weights of loss.

GANs provide a strong model-free, data-driven model for generating realistic scenarios. Based on the required forecast horizon and number of renewable generation plants, We can feed GANs with historical renewable generation data of corresponding dimensions, and then can train the generator and discriminator simultaneously to capture the data distribution of historical observations using \eqref{equ:new_GAN}. Once trained, we can sample from generator to get realistic samples resembling to the samples coming from the joint distribution of renewables generation time-series with pre-defined timing and locational information.

An important benefit of using WGAN is that the training evolution can be continuously evaluated by the output of discriminator, which provides a useful training and evaluation for system operators for real applications. To show that our method preserves this property, we train WGAN on NREL renewable integration dataset and plot the convergence curves of the discriminator's value functions in Fig.~\ref{fig:loss}. The red curves are evaluated on the training set and the orange ones are on the test set. We could observe that the results on the test set can consistently vary with the almost same trend with that of the training set for both wind and solar power, demonstrating that the generative models are well trained. Once the training is completed, we get an optimal generator that can capture the underlying spatial-temporal correlations in renewable generation data. In the next section, we will formally formulate the scenario forecasts problem, where we use the pre-trained GANs for generating a set of future renewables scenarios based on points forecasts.

\section{Scenario Forecasts via GANs}
\label{forecast}
By training GANs to convergence, we get a powerful generator to model the complex spatial and temporal relations  existed in multiple-site renewable generations. However, one more challenging problem is to generate multiple feasible future scenarios. We will show by formulating the scenario forecasts as an optimization problem, we could embed the fitted GANs with any given point forecast method to generate scenarios without dimension and number limitations.

Since the single-site scenario forecasts is a special case of multiple sites, here we give the general problem formulation for multiple sites' scenario forecasts. For a typical multiple renewable power generation sites, assume at timestep $t$, we have some forecasting method to obtain the point forecasts $\hat{\bf p}_{i,j}$ for each power generation site $i$ and each look-ahead time $j$, $i = 1,..., K$, $j=1,..., T$. The point forecasts are given as

\[
\hat{\bf p}_{pred}=
\left[ {\begin{array}{cccc}
	\hat{p}_{1,1} & \hat{p}_{1,2} & \cdots & \hat{p}_{1,T}\\
	\hat{p}_{2,1} & \hat{p}_{2,2} & \cdots & \hat{p}_{2,T} \\
	\vdots &\vdots &\ddots &\cdots\\
	\hat{p}_{K,1} & \hat{p}_{K,2} & \cdots &\hat{p}_{K,T}
	\end{array} } \right]
\]
where $K$ denotes the number of sites and $T$ denotes the forecasting horizon. In this paper, we focus on the scenario forecasting problem, so it is flexible to work with any point forecasts method, e.g., ARIMA model or information from numerical weather prediction (NWP). Given some input noise $Z$, the pre-trained $G(Z)$ generates a possible realization without regarding to the forecast information $\hat{\bf p}_{pred}$ . Based on the generator and the point forecasts, we are interested in forecasting a group of $N$ scenarios
$S =\{S_1,...,S_N\},\; S_i\in \mathbb{R}^{K\times T}$ to represent the uncertainty of renewable generation. Meanwhile, $S$ shall accurately reflect the temporal and spatial dynamics of future generation.

Since the point forecasts  $\hat{\bf p}_{pred}$ only provide a deterministic information for future renewable generations, we take the notion of \emph{prediction interval} to indicate the region around $\hat{\bf p}_{pred}$  that generated scenarios should lie
in~\cite{chen2018unsupervised, wan2014probabilistic}. We describe this interval with an upper bound $U_\alpha(\hat{\bf p}_{pred})$ and a lower bound $L_\alpha(\hat{\bf p}_{pred})$ :
\begin{equation}
L_\alpha(\hat{\bf p}_{pred})=\frac{1}{\alpha} \hat{\bf p}_{pred}; \quad U_\alpha(\hat{\bf p}_{pred})=\alpha \hat{\bf p}_{pred}
\end{equation}
where the hyperparameter $\alpha$ can be interpreted as the prediction confidence or prediction interval.

Since the forecasting scenarios should reflect the forecast information around the point forecast $\hat{\bf p}_{pred}$ , we can first obtain a good starting point for $Z$ by solving the following problem:

\begin{subequations}
	\label{problem:initial}
	\begin{align}
	\min_Z \quad& ||\mathbb{P}_{pred}(G(Z))-\hat{\bf p}_{pred}||_2 \tag{\ref{problem:initial}}\\
	s.t. \quad  & Z\in \mathscr{Z} \label{problem:a_init}\\
	& L_{\alpha}(\mathbf{\hat p}_{pred}) \leq \mathbb{P}_{init}(G(Z)) \leq U_{\alpha}(\mathbf{\hat p}_{pred}) \label{problem:b_init}
	\end{align}
\end{subequations}
where $\hat{\bf p}_{init}$ is randomly sampled from an initial fluctuation interval
$[L_\alpha(\hat{\bf p}_{pred}),\; U_\alpha(\hat{\bf p}_{pred})]$.
Note that our goal is to forecast scenarios that not only can represent the uncertainty of future time, but also can generate realistic time-series that can capture the intrinsic patterns of renewable energy sources at different prediction horizons. According to the loss defined in \eqref{equ:minmax}, larger discriminator output $D(\cdot)$ indicates more realistic samples. To ensure the generated scenarios are realistic with pre-trained generator, we use the objective $-D(G(Z))$ to enforce the generated samples are realistic enough. Meanwhile, we want to constrain generated scenarios within a pre-determined confidence interval according to actual needs of risk management. Using all of the objectives above and pre-trained model $G, \; D$, the scenario forecasts problem can be formulated as a constrained optimization problem:

\begin{subequations}
	\label{problem:main}
	\begin{align}
	\min_Z \quad& - D(G(Z)) \tag{\ref{problem:main}}\\
	s.t. \quad  & Z\in \mathscr{Z} \label{problem:a}\\
	& L_{\alpha}(\mathbf{\hat p}_{pred}) \leq \mathbb{P}_{pred}(G(Z)) \leq U_{\alpha}(\mathbf{\hat p}_{pred}) \label{problem:b}
	\end{align}
\end{subequations}

\begin{figure*}[!t]
	\centering
	\includegraphics[width=1.0 \columnwidth]{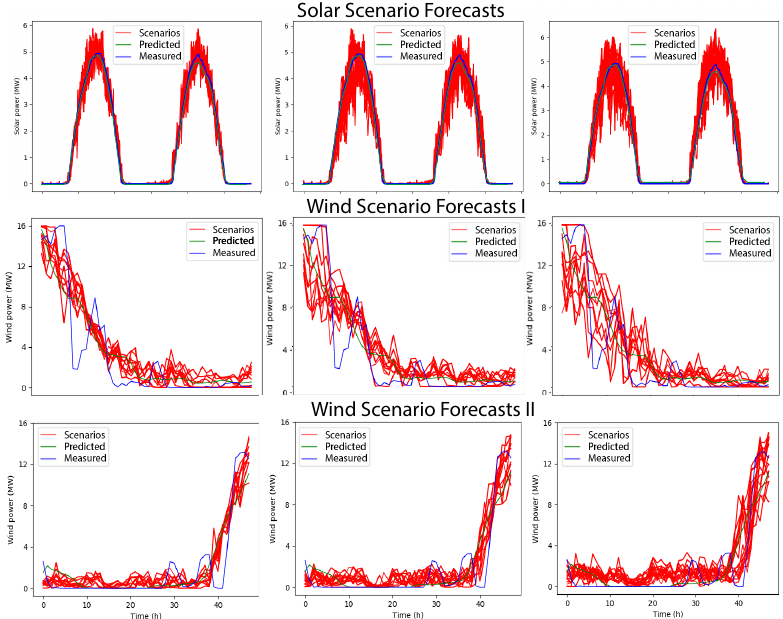}
	\caption{Examples of one solar and two wind power forecasted scenarios with varying PIs of 1.5, 2 and 3 respectively. Note that in first wind example, the point forecasts are not accurate, while increasing PIs could generate scenarios that cover the uncertainties.}
	\label{fig:samples}
\end{figure*}

In order that we can always obtain a good initial $Z$, we set the prediction interval in \eqref{problem:initial} to be slightly smaller than in main objective function \eqref{problem:main}. Since both of the objective and constraints in \eqref{problem:main} are nonconvex, to deal with the inequality constraints, we propose to substitute it into the main objective with two log barriers. Meanwhile, by using the fact that there are multiple local optima in deep neural networks, we can find different initial points $Z_i \in Z$ by solving \eqref{problem:initial} for multiple times,
and get distinct forecasting scenario $\mathbb{P}_{pred}(G(Z^*_i))$ by solving \eqref{problem:main} with log barriers. As the training loss defined in \eqref{equ:new_GAN} incurs $G$ to generate diverse scenarios given different $Z_i$, we are able to obtain a group of distinct yet realistic scenarios that not only can reflect the point forecast information, but also can represent different uncertainty levels according to the actual needs of risk management $\alpha$.

\section{Network Structure and Training Details}
\label{network}
In this section, we describe the network structure that has been used for generating multiple-site short-term scenario forecasts. Note that our structure is flexible in generating scenarios for varying time horizons and number of power plants, it could be served as a plug-in module for generating realistic future scenarios.

Both networks in our model utilize the power of deep learning to represent the temporal and geographical dependencies of scenarios. The generator network starts with fully connected multilayer perceptron and 3 de-convolutional layers to up-sample the input samples $Z$ to generate renewable time-series. The discriminator network has a similar yet reversed structure to distinguish historical samples from generated samples with a scalar output. Sigmoid function is used as the activation function to limit the output of discriminator lying  in the interval $[0,1]$. ReLU and Leaky-ReLU activations are respectively used in the hidden layers of the generator and the discriminator. Dropout is only applied in the output of each hidden layer of the discriminator. Batch normalization can be used to help stabilize training in both the generator and the discriminator, but it changes the form of the discriminator's problem from mapping a single input to a single output to mapping from an entire batch of inputs to a batch of outputs. Since the improved loss objective in \eqref{equ:minmax} is no longer valid in this setting, we can omit or replace the batch normalization by layer normalization in our model structure~\cite{ba2016layer}. All our experiments for scenario forecasts are implemented in Python using open-source machine learning package TensorFlow~\cite{abadi2016tensorflow}.

\section{Computational Experiments}
\label{experiment}
In this section, we describe our experiments and results on two renewable datasets: wind and solar power generations respectively. By first traininig the proposed GANs model on historical datasets, we fix the model weights and implement our scenario forecasts algorithms. We show that the proposed method can forecast short-term scenarios for both single site and spatial-correlated multiple sites. We validate the forecasted scenarios by examining their statistical properties. 

\subsection{Simulation Setup}
In order to test the performance of our proposed framework for scenario forecasts, we build training and validation dataset using power generation data from National Renewable Energy Laboratory~(NREL) Wind and Solar Integration Datasets~\cite{draxl2015wind}. 
Historical power measurements have a resolution of $5$ minutes. We choose $24$ wind farms and $32$ solar power plants located in the State of Washington to use as the datasets for our data-driven method. Historical samples are split into training set and standalone validation set. In general, we randomly select $80\%$ of the data as training set. we also collect the corresponding 24-hour ahead forecast data, which is served as historical observations which are later used for forecasting scenarios based on pre-trained GANs. All renewable power measurements and forecasts are rescaled to $[0, 1]$.

Our scenario generation model is repeatedly fed with the historical samples until the discriminator loss to converge. We keep the training until about $16,000$ iterations to demonstrate the training procedure is stable. The training convergence curves on solar and wind data sets for GANs are shown in Fig.~\ref{fig:loss}. Output of $L_D$ is growing initially since the discriminator could learn to distinguish samples generated by the generator from the real historical samples. 
The generator gradually learns various patterns in historical renewable data. 
After $6,000$ iterations of training, the loss function already converges to near $0$. As the training tends to converge, The generator is able to generate plausible power trajectories with a small $L_D$ and the discriminator can hardly distinguish between generated time-series and real ones. Eventually, the output power scenarios of the generator can represent the stochastic processes of renewable power generation.

\begin{figure}[h]
	\centering
	\includegraphics[width=0.8 \columnwidth]{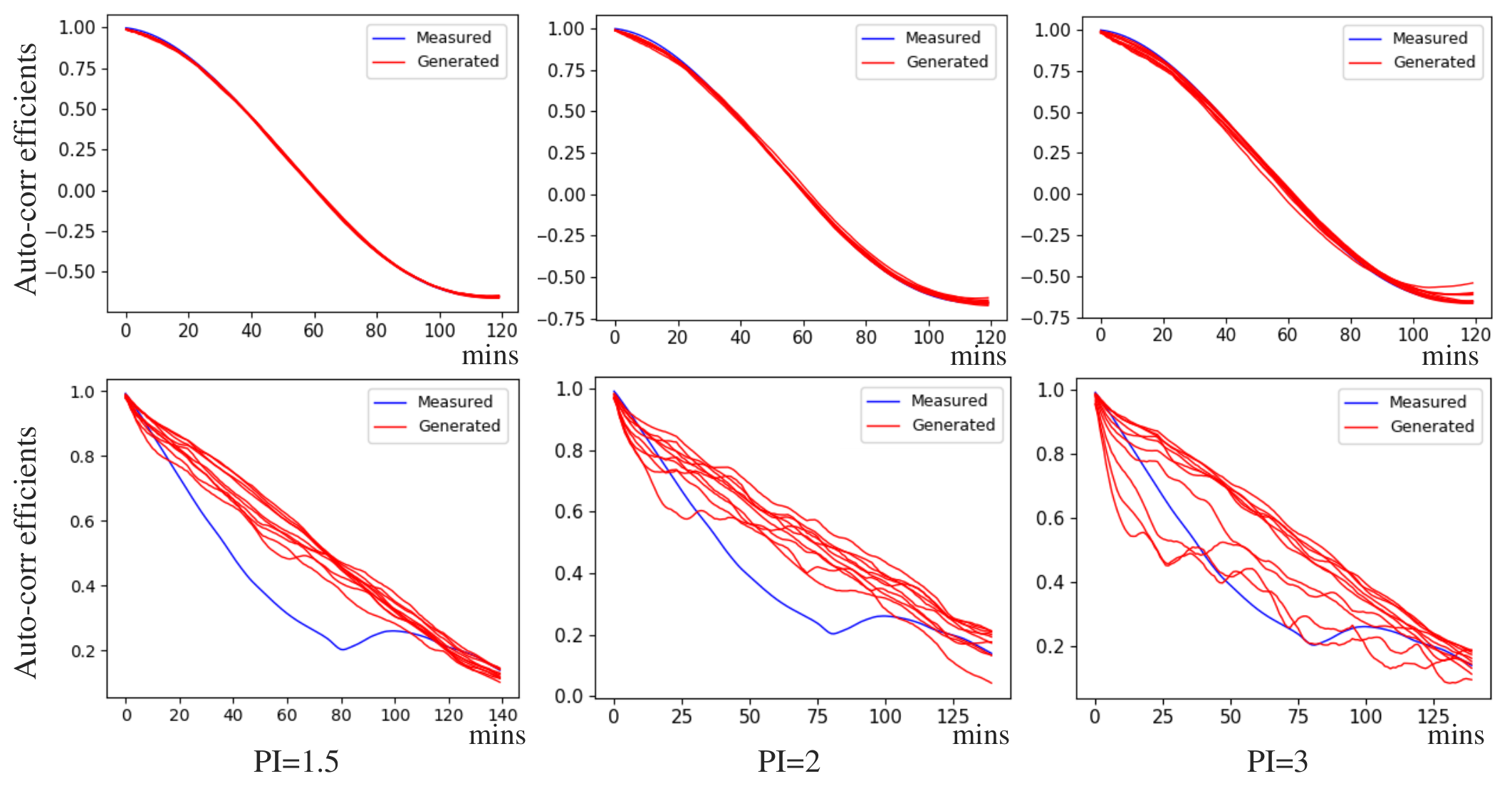}
	\caption{Autocorrelation plots for both measured values and generated scenarios: (top)solar power and (bottom)wind power. The autocorrelation plots correspond to the first two examples shown in Fig.~\ref{fig:samples}.}
	\label{fig:autocorr}
\end{figure}

\subsection{Scenario Forecasts}
For different scenario forecasting tasks, we can use the same GANs model. The framework for using GANs for forecasting scenarios is illustrated in Fig.~\ref{fig:schematic}. The proposed method contains two steps. We can use the first step to model the uncertainty and capture the data distribution of renewable resources. while in the second step, by solving the optimization problem, we can generate a large number of forecasted scenarios.

In this subsection, we validate the proposed method that can forecasting scenarios for a single site of renewable resources. Historical data in geographical proximity is collected as input samples to represent the stochastic generation dynamics. We first show that the proposed method can generate scenarios of different levels of uncertainty for solar and wind power. The size of samples from the training set is composed of two-day data.  The forecasted trajectories with varying PIs of 1.5, 2 and 3 are shown in Fig.~\ref{fig:samples} for both solar and wind cases with 3 randomly selected 2-day samples. By visually inspection, we can observe that the samples generated by proposed methods can correctly capture the hallmark features (e.g., large peak values, daily variations, and ramp events of large fluctuations) of the solar and wind power profiles from the predicted data. By selecting different prediction interval $\alpha$, the forecasted scenarios can represent different degrees of uncertainty in renewable power generation. When the interval level is $\alpha=1.5$, generated scenarios are close to point forecasts, yet fail to cover the realizations; while when $\alpha=3$, generated trajectories could cover the actual power production values, but are less concentrated. 
We note in the first example of wind scenario forecasts, even when the point forecasts are not accurate, a larger prediction interval will help generate diversified scenarios that cover the forecast uncertainties. The prediction interval can be properly tuned based on the real-world applications. Meanwhile, we can further adjust the hyperparameters of the upper and lower boundaries of the log-barrier in optimization problem to generate trajectories that is more in line with the actual needs. 


In order to further verify the generated scenarios' temporal statistics, we calculate and compare samples' autocorrelation. The autocorrelation measures the degree of correlation of a time series between two different periods.  Autocorrelation represents the temporal correlation at a renewable resource. The autocorrelation coefficient $R(h)$ for a given time-series is calculated as
\begin{equation}
\label{equ:autocorr}
R(h)=\sum_{i=1}^{n-h}\frac{(S_i-\mu)(S_{i+h}-\mu)}{\sum_{i=1}^{n}(S_i-\mu)^2}
\end{equation}
where $h$ is the look-ahead time and $S$ represents generated samples or realizations with mean $\mu$.

The temporal correlation of the first two examples shown in Fig.~\ref{fig:samples} are calculated by~\eqref{equ:autocorr}, and the results are shown in Fig~\ref{fig:autocorr}. Scenarios' autocorrelation plots cover the measurements' autocorrelation, indicating the generated scenarios are able to represent the temporal dependence of real time-series. The discussions on prediction intervals also hold similarly on the autocorrelation plots, as the increase of PIs leads to diversified scenarios that have various representation of autocorrelations.

\subsection{Spatial Correlation}
\begin{figure}[h]
	\centering
	\includegraphics[width=0.9 \columnwidth]{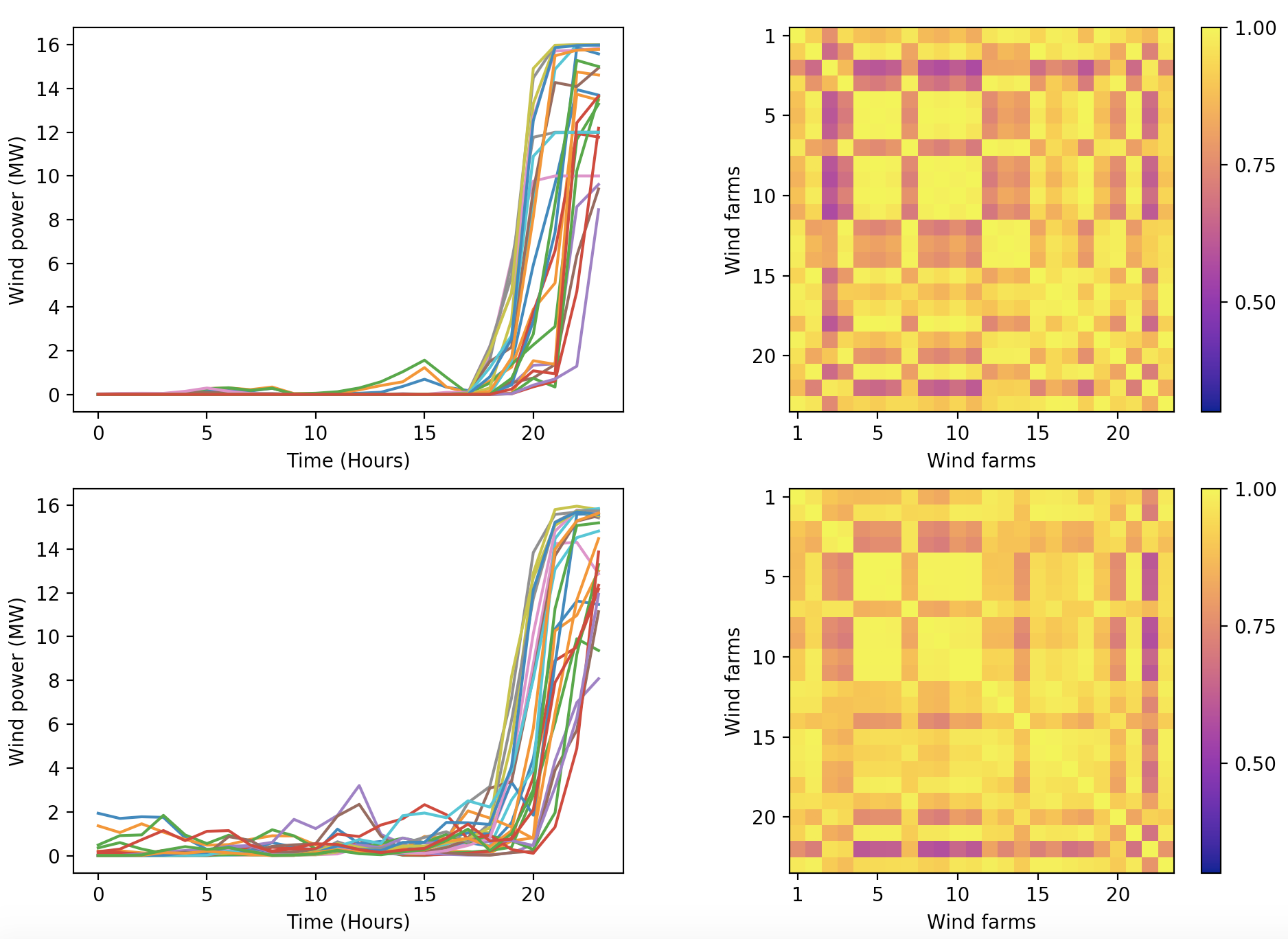}
	\caption{Wind power scenarios and spatial correlation coefficient colormaps for multiple sites: (top) historical data; (bottom) sample generated by our method.}
	\label{fig:spatial_corr}
\end{figure}

For the scenario forecasts of multiple sites, instead of feeding historical data $X$ for a single site to GANs,  we input the model with a real data matrix of size $K\times T$, where $K$ denotes the total number of generation sites, while $T$ denotes the total number of timesteps for each scenario. Here we choose $K = 24, T = 24$ with a resolution of $1$ hour. In order to further examine the correlations between individual locations, we calculate the correlations of the simulated time series and compare the values with those of measured time series. Each row of the correlation matrix shows the correlation between that site (e.g. row 1 represents Wind farm 1) and the other sites, so that the diagonal is composed of ones (the site auto-correlation) and the other terms are the cross-correlation between sites. A sample of 24 wind farms' real power generations  and forecasting scenarios along with the correlation matrix between different sites are plotted in Fig. \ref{fig:spatial_corr}. By visual inspection we find that their dynamic behaviors are similar to each other. The spatial and temporal correlations in the real data (again, not seen in the training stage) are correctly preserved by our forecasted scenarios. From the spatial correlation coefficient colormaps of these two group of data, we can see the generated scenarios preserve the relative correlation between any pair of wind farms. Such results indicate that our proposed method could preserve both the spatial and temporal correlations, and can be used as a plug-and-play module for multi-site planning and operation problem.

\begin{figure}[h]
	\centering
	\includegraphics[width=1.0 \columnwidth]{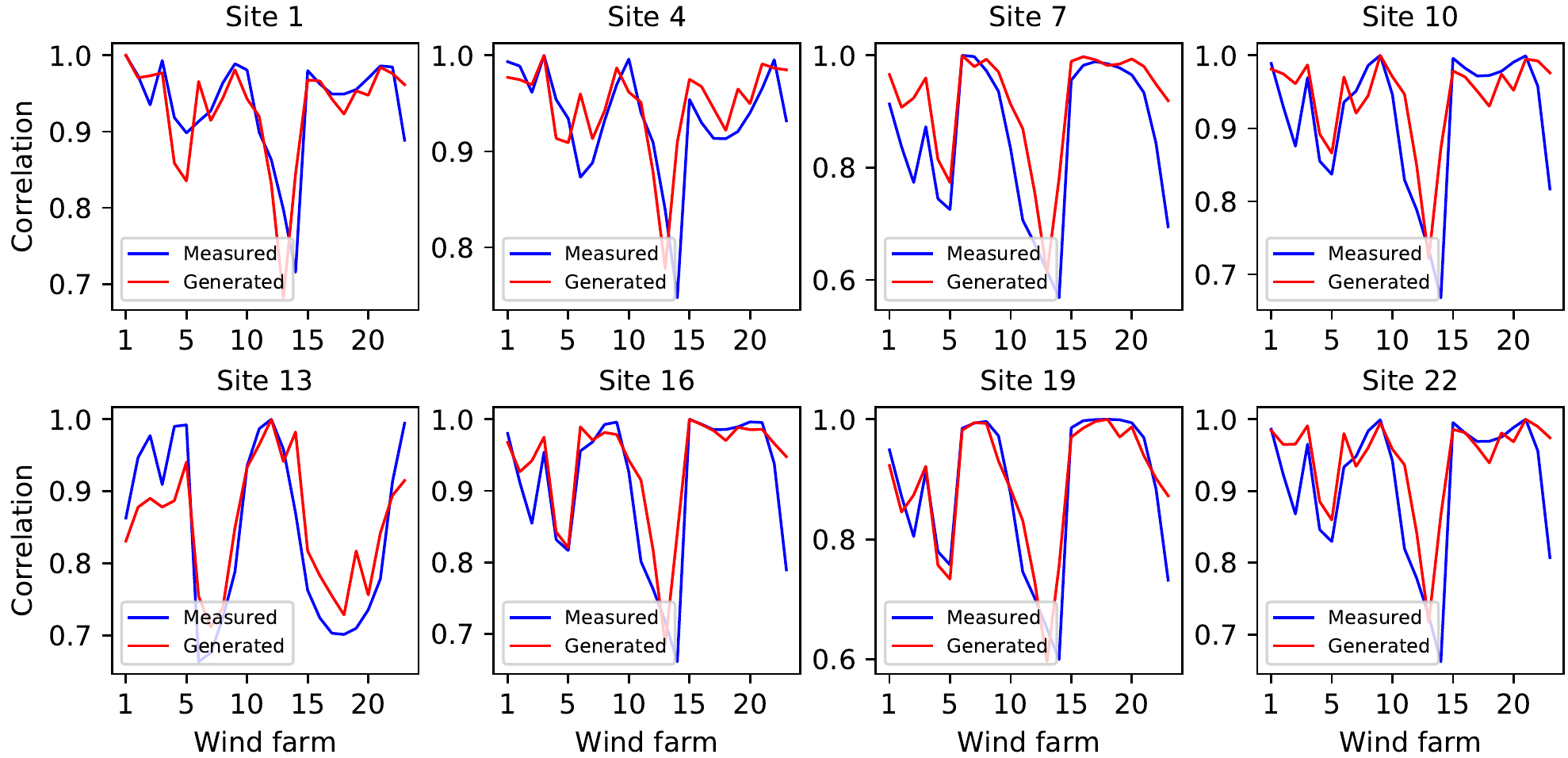}
	\caption{Correlation between a given set and other sites, with a comparison between real data and generated scenarios.}
	\label{fig:correlation_example}
\end{figure}

For purposes of illustration, some elements of the correlation matrix for a few sites are also shown in Fig.~\ref{fig:correlation_example} for both the measured time-series and the forecasted scenarios. Each pair of curves basically maintains a consistent trend compared to every other site. The results show that the generated scenarios using proposed method agree with the ground truth, showing that spatial correlations between different sites are correctly retained.

\begin{figure}[h]
	\centering
	\includegraphics[width=0.8 \columnwidth]{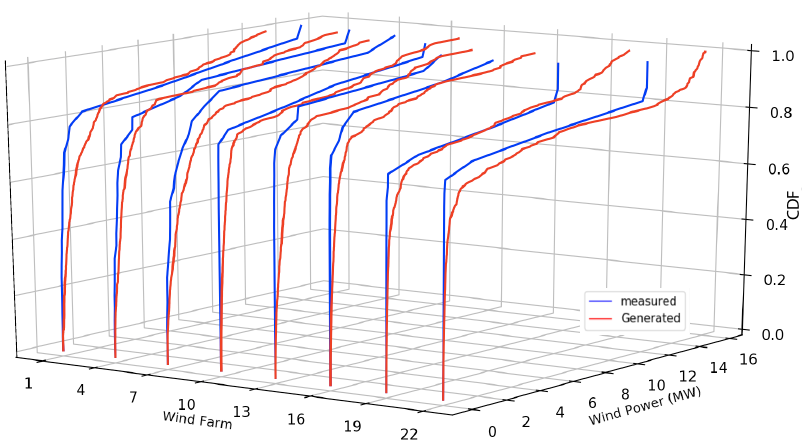}
	\caption{Cumulative distribution function (CDF) of measured wind power versus CDF of generated scenarios via our method.}
	\label{fig:cdf}
\end{figure}

We further verify that generated time-series have the same statistical properties as the measured data. We generate 50 scenarios for a group of 24 wind farms. We compare the CDF  of the true realizations and the generated scenarios in Fig.~\ref{fig:cdf}. For the sake of simplicity, we only select some of them for display. It is clear the methodology for different sites has the capability to generate samples with the correct marginal distributions that are basically the same as the predicted scenarios. The power generation and fluctuations at these sites are  consistent with the measured data. The level of different generation capacities and the magnitude of different fluctuations can be correctly captured.

\section{Conclusions}
\label{conclusion}
This paper proposes a novel method to forecast scenarios for renewables power generation processes based on deep generative models. Capable of working with any off-the-shelf point forecast methods, the proposed algorithm not only characterizes the uncertainty associated with renewable energy resources for both single site and spatially correlated multiple-site cases, but also generates realistic short-term scenarios without any number or forecast horizon restrictions.  Comprehensive simulations carried out for different case studies show the effectiveness of the proposed methodology. With high reliability and high flexibility, the proposed approach can be used to directly generate a large number of scenarios. Various statistical methods validate the superiority of proposed approach: the marginal distribution associated with each renewable power stochastic process is retained by the generated scenarios; the temporal correlations are characterized by autocorrelations at each renewable stochastic process; the spatial correlations are verified by cross-correlations among different sites. Our method can serve as a promising pipeline for generating high quality scenarios that reflect the intrinsic patterns and data distribution of renewablge generations.

\bibliographystyle{IEEEtran}
\bibliography{bib}

\end{document}